\documentclass[oneside]{amsart}

\usepackage[T1]{fontenc}

\usepackage[latin1]{inputenc}

\usepackage{amsmath,amsthm,amsfonts,amssymb,latexsym, mathrsfs}

\makeatletter

\theoremstyle{plain}
\newtheorem{thm}{Theorem}[section]

\numberwithin{equation}{section} 

\numberwithin{figure}{section} 

\theoremstyle{plain}
\newtheorem{cor}[thm]{Corollary} 

\theoremstyle{definition}
\newtheorem{defn}[thm]{Definition}

\theoremstyle{plain}
\newtheorem{lem}[thm]{Lemma} 

\theoremstyle{plain}
\newtheorem{prop}[thm]{Proposition} 

\theoremstyle{plain}
\newtheorem{fact}[thm]{Fact}

\theoremstyle{plain}

\newtheorem{rem}[thm]{Remark}

\newtheorem{ques}[thm]{Question}

\makeatletter

\usepackage{amssymb}

\usepackage{amsmath}

\newcommand{\N}{\mathbb{N}}
\newcommand{\Z}{\mathbb{Z}}
\newcommand{\Q}{\mathbb{Q}}

\def \<{\langle}
\def \>{\rangle}

\usepackage{amsfonts}

\usepackage{amscd}

\usepackage{latexsym}

\usepackage{epsfig}

\makeatother

\def\Ind#1#2{#1\setbox0=\hbox{$#1x$}\kern\wd0\hbox to 0pt{\hss$#1\mid$\hss}
\lower.9\ht0\hbox to 0pt{\hss$#1\smile$\hss}\kern\wd0}

\def\Notind#1#2{#1\setbox0=\hbox{$#1x$}\kern\wd0\hbox to 0pt{\mathchardef
\nn=12854\hss$#1\nn$\kern1.4\wd0\hss}\hbox to
0pt{\hss$#1\mid$\hss}\lower.9\ht0 \hbox to
0pt{\hss$#1\smile$\hss}\kern\wd0}

\def\thind{\mathop{\mathpalette\Ind{}}^{\text{\th}} }

\def\nthind{\mathop{\mathpalette\Notind{}}^{\text{\th}} }

\def\uth{\text{U}^{\text{\th}} }

\newcommand {\thorn} {\text{\th}}

\newcommand {\tr}[1] {\thorn\text{-rank}\big(#1\big)}

\newcommand {\onto} {\twoheadrightarrow}

\newcommand{\eq}{^{\operatorname{eq}}}

\newcommand{\tp}{\operatorname{tp}}

\newcommand{\dcl}{\operatorname{dcl}}

\newcommand{\scl}{\operatorname{scl}}

\newcommand{\abs}[1]{\lvert#1\rvert}

\begin{document}

\title{Thorn independence in the field of real numbers with a small multiplicative group}

\author{Alexander Berenstein}

\address{Universidad Nacional de Colombia,
Cra 30 No 45-03,
Bogot\'{a}, Colombia}
\email{ajberensteino@unal.edu.co}
\urladdr{www.matematicas.unal.edu.co/\textasciitilde aberenst}

\author{Clifton Ealy}

\address{University of Illinois at Urbana-Champaign,
1409 West Green Street,
Urbana, Illinois 61801}

\email{clif@math.uiuc.edu}

\urladdr{www.math.uiuc.edu/\textasciitilde clif}

\author{Ayhan G\"unayd\i n}

\address{University of Illinois at Urbana-Champaign,
1409 West Green Street,
Urbana, Illinois 61801}

\email{gunaydin@math.uiuc.edu}

\urladdr{www.math.uiuc.edu/\textasciitilde gunaydin}



\date{June 1st, 2007}

\thanks{The authors would like to thank Lou van den Dries for his close reading and helpful comments.}

\begin{abstract}

We characterize \th-independence in a variety of structures, focusing on the field of real numbers expanded by predicate defining a dense multiplicative subgroup, $G$, satisfying the Mann property and whose $p$th powers are of finite index in $G$. We also show such structures are super-rosy and eliminate imaginaries up to codes for small sets.

\end{abstract}

\pagestyle{plain}

\maketitle

\section{Introduction}

We build on results of van den Dries and G\"unayd\i n in \cite{Mann property}. There the authors investigate the model theory of pairs $(K,G)$ where $K$ is either an algebraically closed field or a real closed field, and $G$ is a multiplicative subgroup of $K^\times$ with the Mann Property. 
While the definition of the Mann property is somewhat lengthy (and we postpone the precise definition to Section 5), roughly the Mann Property is a condition insuring that linear equations have few solutions in $G$.  Among other things, the Mann property implies that $G$ is \emph{small} (in a technical sense defined below).  Moreover, such groups are quite natural.  Any group contained in the divisible hull of a finitely generated group, i.e. any finite rank group, has the Mann property.  

In the case where $K$ is real closed (henceforth we distinguish this case by referring to $K$ as $R$), the additional hypothesis that $G$ is  a dense subgroup of $R^{>0}$ is used. 

Among other results, van den Dries and G\"unayd\i n obtain good descriptions of the definable sets in both cases and a good description of dimension when $K$ is algebraically closed, assuming $G$ is $\omega$-stable.  In particular, the pair $(K,G)$ is shown to be $\omega$-stable of Morley rank $\omega$.  

We extend the results of \cite{Mann property} by obtaining a description of dimension for $R$ real closed and 
$G$ such that for each prime number, $p$, the subgroup of $G$ consisting of $p$th-powers has finite index in $G$.  To do this, we need to refine slightly the description of definable sets, focusing on a certain collection of definable sets we call ``basic small'', and introduce the notion of \th-rank.  In particular, we prove that the pair $(R,G)$ is super-rosy of \th-rank $\omega$.  We then use this fact to obtain some partial results about elimination of imaginaries.

Now we state these results precisely.

\begin{thm}\label{Main Theorem 1}
Let $R$ be a real closed field and $G$ a dense  subgroup of $R^{>0}$ with the Mann property and such that for each prime number, $p$, the subgroup of G consisting of $p$th-powers in $G$ has finite index in $G$. Then in the language of ordered rings augmented with a unary predicate for $G$, we have

(1) $G$ has \th-rank $1$, and

(2) $(R,G)$ has \th-rank $\omega$.

\noindent Hence, $(R,G)$ is super-rosy.
\end{thm}

\begin{thm}\label{Main Theorem 2}
Let $(R,G)$ be as in the previous theorem. Enlarge $(R,G)$ by adding sufficiently many sorts of $(R,G)^{\eq}$ so that the resulting structure has a code for every basic small subset of $R^k$, for each $k$. Then this structure eliminates imaginaries.
\end{thm}

While our primary interest is in subgroups of $\mathbb R$ with the Mann property, we obtain Theorems \ref{Main Theorem 1} and \ref{Main Theorem 2} as applications of a more general result: 

\begin{thm} \label{Main Theorem 3}
Suppose that  $(R,+,\dots)$ is an o-minimal expansion of a group in the language $\mathscr L$.  Consider the expansion $\mathfrak R=(R,\mathcal G,+,\dots)$ in the language $\mathscr{L}_\mathcal G=\mathscr L \cup \{\mathcal G\}$ where $\mathcal G$ is a unary predicate. Suppose that for each $\mathfrak{R}'=(R',\dots)$ with $\mathfrak R'\equiv\mathfrak R$:
\begin{enumerate}
\item $\mathcal G(R')$ is small, and contained in some interval, $(a,\infty)\subseteq R'$, in which it is dense.

\item Each $\mathscr{L}_\mathcal G$-formula $\psi(x)$  is equivalent to a boolean combination\footnote{Throughout the paper, we use ``boolean combination of . . .'' to mean ``an element of the ambient boolean algebra generated by . . .''.} of formulas of the form $\exists \vec y \big(\mathcal G(y_1)\wedge\dots\wedge \mathcal G(y_j)\wedge \varphi(x,\vec y)\big)$ where $\varphi$ is an $\mathscr L$-formula.

\item For each tuple $\vec a$ from $R'$ and $\mathbb D\subseteq \mathcal G(R')^n$, definable over $\vec a$, there are an $\mathscr L$-definable set $\mathbb E$, and a definable $\mathbb S$, which is a dense subset of $\mathcal G(R')^n$, with $\mathbb E$ and $\mathbb S$ over $\vec a$, such that $\mathbb D =\mathbb E\cap \mathbb S$.   Furthermore, when $n=1$, $\mathbb D$ can be written as a finite union of such $\mathbb E \cap \mathbb S$, where $\mathbb S$ is, in addition, $\emptyset$-definable. 
\end{enumerate}
Then  $\mathfrak R$ is super-rosy of \th-rank less than or equal to $\omega$ and \th-rank of $\mathcal G(R)$ is 1.  Moreover, if $\mathfrak R$ includes a field structure, the \th-rank of $\mathfrak R$  equals $\omega$.

\end{thm}

For the definition of {\it small}, see \ref{small}.
\smallskip

The reader will note that if conditions (1) and (2) hold in a given model, they hold in any elementarily equivalent model, and if condition (3) holds in a sufficiently saturated model, it holds in any elementarily equivalent model.  The reader will further note that condition (3) above seems quite technical.  In many cases, a much more natural (and stronger) condition holds.  Namely, 

$(3)'$\emph{ For each definable  $\mathbb D\subseteq \mathcal G(R)^k$ there is an $\mathscr L$-definable set $\mathbb E$ such that $\mathbb D = \mathbb E\cap \mathcal G(R)^k$}.  

However, in cases that are of particular interest to us, such as $\mathfrak R=(\mathbb R,\mathcal G,+,\cdot)$ and $\mathcal G(R)=2^\mathbb Z 3^\mathbb Z$, $(3)'$ fails.  To understand why $(3)$ is not as unnatural as it may first appear, the reader may skip ahead to Section \ref{bounded case}.

\begin{thm} \label{Main Theorem 4}
Let $\mathfrak R$ be as in the previous theorem.  Enlarge $\mathfrak R$ by adding sufficiently many sorts of $\mathfrak R^{\eq}$ so that the resulting structure has a code for every basic small subset of $R^k$. Assume in addition, given any set of parameters $A$, and any interval $I$  defined over $A$, that $\scl(A)\cap I$ is not contained in any small set (see \ref{small} and \ref{small closure} for the appropriate definitions). Then this structure eliminates imaginaries. 
\end{thm}

In addition to applying to structures satisfying the conditions of Theorem \ref{Main Theorem 1},  Theorems \ref{Main Theorem 3} and \ref{Main Theorem 4} also apply to the structures studied in  \cite{Dense Pairs}, namely dense pairs of o-minimal expansions of ordered abelian groups.  Also we note that we answer the question of Miller and Speissegger from \cite{o-minimal open core} of whether $(\mathbb R, 2^Z 3^Z)$ has o-minimal open core.  (See the end of Section 3 for definitions and the main part of the proof, and Section 5 for its application to expansions of the reals by groups with the Mann Property.)

\medskip

\subsection*{Conventions and Notation}

An $\mathscr L$-structure, e.g. $\mathfrak R=(R,+^{\mathfrak R},\cdot^{\mathfrak R},<^{\mathfrak R},0^{\mathfrak R},1^{\mathfrak R})$, consists of an underlying set, e.g. $R$, together with an interpretation of each symbol from the language, e.g $+^{\mathfrak R},\cdot^{\mathfrak R},<^{\mathfrak R}, 0^{\mathfrak R},1^{\mathfrak R}$.  We drop the superscripts when no confusion results. Capital letters in the Fraktur font, e.g. $\mathfrak{M}$ and $\mathfrak{R}$, indicate structures.  The universes of these structures are denoted by the corresponding capital letters in the normal font.  For instance, $M$ and $R$ are the respective universes of the structures above.

We use the letters $x,y,z,w$ as variables, and the letters, $a,b,c,$ etc., to indicate elements of the universe of a structure.  We distinguish between elements from $M$ and tuples from $M^n$  by using vector notation for tuples.  For example, $\vec x, \vec y$ and $\vec a, \vec b$ as opposed to $x,y$ and $a,b$.    

We use $\varphi$, $\psi$, and $\theta$ to indicate formulas.  When no confusion results, we suppress the parameters, writing, for instance, $\varphi(\vec x)$ even when the formula is not over the empty set. Likewise, when we say definable, we mean definable with parameters.

To save ourselves from constantly worrying about the length of our tuples, when  $\vec x$ is an $n$-tuple, we write $M^n$ as $M_{\vec x}$. The set defined by a formula $\varphi(\vec x)$ is denoted by $\varphi(M_{\vec x})$.

We use capital letters in blackboard bold to indicate definable sets, e.g. $\mathbb{D},\mathbb{E}$, with the exceptions of $\mathbb{N}$, $\mathbb Q$, and $\mathbb R$, which are the sets of natural numbers, rational numbers, and real numbers, respectively. We denote the complement of $\mathbb D$ as $\mathbb D^c$. We use $f$, possibly with subscripts, for definable functions. Also $\alpha,\beta$ will always indicate ordinals, $m,n$ will always indicate natural numbers, and $p$ will always indicate a prime number. 

If we wish to emphasize that a definable set is defined with parameters, we write the parameters as a subscript.  For example, suppose $\psi(\vec y)$ defines $\mathbb E$ and $\varphi(\vec x)$ defines $\mathbb D$, where we have suppressed the parameters in both $\psi$ and $\varphi$.  If we then wish to emphasize that $\varphi$ uses a parameter $\vec e\in M_{\vec y}$, we write $\mathbb D_{\vec e}$.  For instance, we write $$\exists \vec y(\psi(\vec y)\land \varphi(M_{\vec x},\vec y))$$ as $$\bigcup_{\vec e \in \mathbb E}\mathbb D_{\vec e}.$$

For a set $C$, we denote by $\mathscr P(C)$ the power set of $C$.

\subsection*{Definitions and Preliminaries}
Now we introduce some definitions that we use in the remainder of the paper, together with some propositions from other papers which we also use.

\begin{defn}
Fix a theory, $T$, and a sufficiently saturated model $\mathfrak M \models T$.  We work in $\mathfrak{M}^{\eq}$. Let $\varphi(\vec x,\vec y)$ be a formula without parameters, let $\vec b\in M_{\vec y}^{\eq}$, and let $C$ be a set of size less than the degree of saturation of $\mathfrak{M}$.

For $k\in \mathbb N$, the formula $\varphi(\vec{x},\vec b)$ is said to \emph{k-\th-divide} over $C$ if there is $D\supseteq C$ such that $\tp(\vec b/D)$ is not algebraic and the set of formulas
$\{\varphi(\vec x,\vec b') : {\vec b'\models\tp(\vec b/D)}\}$ is
$k$-inconsistent. The formula is said
to \emph{\th-divide} over~$C$ if it $k$-\th-divides for some $k$.

The partial type $\pi(\vec x,\vec b)$ is said to \emph{\th-fork} over $C$ if it
implies a disjunction of formulas (with arbitrary parameters), each of which \th-divides over~$C$.
\end{defn}

\medskip

We have defined what it means for a formula to \th-divide over a set $C$.  Sometimes, when the particulars of $C$ are not important, we will simply say that a formula \th-divides.
\begin{rem}
By compactness, if $\varphi$ $k$-\th-divides, there is always a single formula $\theta(\vec y,\vec d)\in \tp(\vec b/D)$ such that the set of formulas
$\{\varphi(\vec x,\vec b') : \mathfrak M\models\theta(\vec b',\vec d)\}$ is
$k$-inconsistent.

Also by compactness, if $\pi(\vec x, \vec b)$ implies a disjunction of formulas that \th-divide, $\pi$ implies a finite disjunction of such formulas.
\end{rem}

\begin{defn}
Let $A,B,C\subset M$ be smaller than the degree of saturation of $\mathfrak M$. Then $\thind$ is defined
as follows: $A\thind_CB$ if and only if $\tp(\vec a/BC)$ does not \th-fork
over~$C$ for any tuple $\vec a$ from $A$.
If $A\thind_CB$ we say that \emph{$A$ is \th-independent from $B$ over $C$.}  If it is clear from context, we will often just say \emph{independent} for \th-independent.
\end{defn}

\begin{defn}
A theory $T$ such that $\thind$ is symmetric for $T$ is called \emph{rosy}.
\end{defn}

Alternatively, rosiness could be defined in terms of local \th-ranks being finite.  However, we will not have need of any local ranks as the situation in which we find ourselves allows for a global \th-rank, as defined below.  

When working with an independence relation, we can define its foundation rank. For \th-independence we have:

\begin{defn}

Let $p(x)\in S(A)$.
For $\alpha$ an ordinal, we define $\uth(p )\geq \alpha$ inductively on $\alpha$.
\begin{enumerate}
\item $\uth(p(x))\geq 0$.
\item If $\alpha=\beta +1$, we define $\uth(p(x))\geq \alpha$ if there is a tuple $a$ and a type  
$q(x,y)$ over $A$ such that $q(x,a)\supset p(x)$, $\uth(q(x,a))\geq \beta$ and $q(x,a)$ \th-forks over $A$.
\item If $\alpha$ is a limit ordinal, then $\uth(p(x))\geq \alpha$ if $\uth(p(x))\geq \beta$ for all $\beta<\alpha$.
\end{enumerate}
\end{defn}

\begin{rem}
It is perhaps worth noting that in a theory that is not rosy, \th-forking may still be symmetric if one restricts the sorts that one considers. If thorn independence satisfies symmetry when restricted to the real sorts, one calls the theory \emph{real-rosy}.  For instance, the theory of algebraically closed valued fields is not a rosy theory, but \th-forking, restricted to the field, residue field, and value group sorts, is an independence relation. Thus ACVF is real-rosy \cite{EO}.  
\end{rem}

\begin{defn}
\th-rank is the least function taking values in $\text{On}\cup\{\infty\}$ satisfying the following:

(1) $\tr{\varphi(\vec x,\vec b)}\geq 0$ if $\varphi(\vec x,\vec b)$ is consistent. 

(2) $\tr{\varphi(\vec x,\vec b)}\geq \alpha+1$ if there is $\psi(\vec x, \vec c)$ that \th-divides over $\vec b$, such that $\psi(\vec x, \vec c) \vdash \varphi(\vec x, \vec b)$ and $\tr{\psi(\vec x,\vec c)}\geq \alpha$.

(3) For $\lambda$ a limit ordinal, $\tr{\varphi(\vec x,\vec b)}\geq \lambda$ if $\tr{\varphi(\vec x,\vec b)}\geq \alpha$ for all $\alpha<\lambda$. 
 
\end{defn}

\noindent The relation between \th-rank and $\uth$-rank is given by the following (\cite{EO}):

\begin{fact}\label{U-rank is smaller}
For any type, $p$, $\uth(p)\leq \text{min}\{\tr{\varphi} : \varphi\in p\}$.
\end{fact}

In analogy with simple and stable theories, we make the following definition (which could be equivalently stated in terms of $\uth$-rank, see \cite{EO}):

\begin{defn}
A complete theory is {\it super-rosy} if every formula has ordinal \th-rank.
\end{defn}

The corollary of the Coordinatization Theorem of \cite{On1} stated below will simplify our proof of super-rosiness:

\begin{cor}
\label{cor-to-coord} Given a complete theory $T$, if every formula in one free variable $\varphi(x,\vec b)$ has ordinal \th-rank, then $T$ is super-rosy. 
\end{cor}

\begin{defn}\label{small}
Let $\mathfrak{M}:=(M,\ldots)$ be an ordered structure. A definable set $\mathbb D\subset M^k$ is
{\it large} iff there is some $m$, an interval $I\subseteq M$ and a function $f:\mathbb D^m\onto I$.

A definable set $\mathbb S$ is {\it small} iff it is not large.
\end{defn}

Note that this definition of small differs from the conventions of \cite{Mann property}. There the adjective ``small'' also applies to sets that are not definable, but does not apply to subsets of $M^n$ for $n>1$.  In addition, in \cite{Mann property}, the notion of small set is defined for arbitrary, possibly unordered, structures.   One of the cases we wish to consider, however, is dense pairs of ordered abelian groups.  In this setting, a bounded interval would be small under the definition of \cite{Mann property}.  Our definition for small, when restricted to definable subsets of a model $(R,G)$ satisfying the hypotheses of Theorem \ref{Main Theorem 3} will turn out to be $G$-small, as defined in \cite{Dense Pairs}.  When $R$ in addition has a field structure all three definitions will coincide (for definable subsets of $R$).

\begin{fact}\label{fact 1}
Let $\mathfrak{M}$ be an o-minimal structure.  Let $\{\varphi(M,\vec{a})\}_{\vec a \in \mathbb{A}}$ be a definable family of subsets of $M$, each of which by o-minimality may be decomposed into a finite union of points and open intervals.  Then the minimal number of points and the minimal number of open intervals in any such decomposition are definable properties of $\vec{a}$.
\end{fact}

Unless stated otherwise, $\mathscr L$ denotes a language extending the language of ordered abelian groups, $\mathcal G$ a unary predicate not in $\mathscr L$, $\mathfrak R=(R,G)$  denotes a structure satisfying the conditions of Theorem \ref {Main Theorem 3}, although one may think of $\mathfrak R$ as a structure satisfying the conditions of Theorem \ref{Main Theorem 1}.  Following our normal conventions, we should refer to the set defined by $\mathcal G(x)$ as $\mathbb G$, but we simply write it as $G$.  We use $\mathfrak{R}|_{\mathscr L}$ to denote the reduct of $\mathfrak R$ to $\mathscr L$. 

\section{Small Sets}

We first make a definition and a technical observation.

\begin{defn}
A  \emph{k-valued function}, $F:A\stackrel{k}{\longrightarrow}B$ is a function from $A$ to $\{S \in \mathscr P(B):|S|\leq k\}$. The {\it graph} of such an $F$ is $\{(a,b)\in A\times B:b\in F(a)\}$, and its {\it image} is $\{b\in B:b\in F(a)\text{ for some $a\in A$}\}$.  If $F:\mathbb D\to\mathbb E$ where $\mathbb D\subseteq R^m,\mathbb E\subseteq R^n$, then we say $F$ is {\it definable in} $\mathfrak R$, if its graph is. 
\end{defn}

We define the composition of such functions as follows:

\begin{defn}
Consider $F_1:A\stackrel{k_1}{\longrightarrow}B$ and $F_2:B\stackrel{k_2}{\longrightarrow}C$.  We define $F_2\circ F_1:A\stackrel{k_3}{\longrightarrow}C$ by setting $F_2\circ F_1(a):=\{c:\exists b \in F_1(a) \textrm{ and }c\in F_2(b)\}$, where $k_3:=k_1\cdot k_2$.
\end{defn}

\begin{lem} \label{many to one}
Let $\mathfrak{M}=(M,<,\ldots)$ be any ordered structure, $\mathbb E,\mathbb F$ be
definable subsets of $M^m$, $M^n$ respectively, and $F:\mathbb E\stackrel{k}{\longrightarrow} \mathbb F$ be a $k$-valued function. Then there is a function $f:\mathbb E^k\to\mathbb F$ with the same image as $F$. If there are two definable elements of $\mathbb E$ then $f$ has the same
parameters as $F$.
\end{lem}

\begin{proof}
Pick distinct $a_1, a_2$ definable elements contained in $\mathbb E$ (adding parameters if necessary). Suppose
that $e\in \mathbb E$ is not equal to $a_1$. Set $f((e,a_1,\dots,a_1))$ to be
the least element of $F(e)$, set $f((a_1,e,a_1,\dots,a_1))$ to be
the second least element of $F(e)$, etc. Now suppose that $e=a_1$. 
Set $f((e,a_2,\dots,a_2))$ to be the least element of $F(e)$, etc.
Finally, for any $\vec e\in\mathbb E^k$ on which $f$ is not yet defined,
set $f(\vec e)$ equal to the least element of $F(a_1)$.
\end{proof}

Let us make a couple of observations about the notion of small as it applies in the setting of groups.  Let $(M, +, \dots)$ be an expansion of a group.  Then the complement of any small set, $\mathbb S$, is large. This can be seen, for instance, by considering the map $f:M^2\to M$ given by $(m_1,m_2)\mapsto m_1+m_2$.  Suppose some element, $m_0\in M$, is not in the image of $(\mathbb S^c)^2$ under $f$.  Then 
$$m_0\in\bigcap_{m\notin\mathbb S}\mathbb S+m.$$ 
Thus, $m_0-\mathbb S$ contains $\mathbb S^c$. 
Now the $2$-valued function
$\mathbb S\stackrel{2}{\longrightarrow} M, \quad s\mapsto\{s,m_0-s\}$
witnesses that $\mathbb S$ is large, which is a contradiction. 
Actually we need a stronger statement:

\begin{lem}\label{small in group}
Let $(M,+,<,\dots)$ be an expansion of an ordered group, and $I=(a,b)\subseteq M$ be a nonempty interval, and $\mathbb S\subseteq M$ a small set. Then $I\setminus\mathbb S$ is large. 
\end{lem}

\begin{proof}
 Let $f:M^2\to M$ be defined as in the previous paragraph. Let $J=(a+b,2b)$. We show that $f\big((I\setminus\mathbb S)^2\big)\supseteq J$.  For a contradiction, let $m_0\in J\setminus f\big((I\setminus\mathbb S)^2\big)$. Then, reasoning as above, $-(\mathbb S\cup I^c)+m_0\supseteq I\setminus\mathbb S$. Noting that $$I^c+m_0=(-\infty,-b+m_0)\cup( -a+m_0, \infty),$$ we see that this yields $-\mathbb S+m_0\supseteq(-b+m_0,b)$, contradicting the smallness of $\mathbb S$.
\end{proof}

\begin{defn}
We say a definable set $\mathbb D$ is {\it small in an interval} $I$ if $\mathbb D\cap I$ is small.
We say a definable set $\mathbb D$ is {\it cosmall in an interval} $I$ if $\mathbb D^c\cap I$ is small. 
\end{defn}

\noindent
Here we return from considering arbitrary ordered groups to the setting of Theorem \ref{Main Theorem 3}.

\begin{defn}
A definable set $\mathbb X$ is {\it basic} if it is defined by a formula of the form $\exists \vec y(\mathcal G(\vec y)\wedge\varphi(\vec x,\vec y))$ where $\varphi(\vec x,\vec y)$ is a formula in $\mathscr L$, and by $\mathcal G(\vec y)$, we mean $\mathcal G(y_1)\land\cdots\land \mathcal G(y_n)$.  Furthermore, we will refer to formulas of the form $\exists \vec y(\mathcal G(\vec y)\wedge\varphi(\vec x,\vec y))$ as \emph{basic formulas}.
\end{defn}

\begin{rem}
Note that a set is basic if and only if it can be written as $$\bigcup_{\vec g\in G^n}\varphi(R_{\vec x},\vec g).$$  where $\varphi$ is an $\mathscr L$-formula. Note also that finite unions and intersections of basic sets are again basic.  In particular, an interval intersect a basic set is again a basic set.
\end{rem}


\medskip

For our purposes the above characterization of definable sets is not quite sufficient; we obtain a more detailed description in the case of definable subsets of $R$ (as opposed to $R^n$).

First we need to prove that if $f_1$ and $f_2$ are functions $R^n\to R$ definable in $\mathscr L$ then $\bigcup_{\vec g\in G^n}(f_1(\vec g),f_2(\vec g))$ is a finite union of intervals.  This is clear when $f_1$ and $f_2$ are functions in one variable.  In general, it is slightly less clear. However, it is a consequence of the cell decomposition theorem for o-minimal structures and the following two lemmas.

The first of the two lemmas shows that subsets of $G^k$ are in a sense well approximated by $\mathscr L$-definable sets.   We already know that for any such set, $\mathbb D$, there is an $\mathscr L$-definable set $\mathbb E$ such that $\mathbb D$ is dense in $\mathbb E\cap G^k$.   It is not the case that $\mathbb D$ will necessarily be dense in $\mathbb E$.  For instance, let $(R,G):=(\mathbb R,2^\mathbb Q)$.  Consider the plane, $\mathbb P\subset \mathbb R^3$ defined by $z-3y=0$.  Let $\mathbb D:=\mathbb P\cap G^3$.  Then $\mathbb D$ is just the copy of $G$ lying on the $x$-axis, and not dense in $\mathbb P$.   Clearly, in this example, had we chosen $\mathbb E$ as the $x$-axis, rather than the plane $\mathbb P$ we would have obtained the density we desired.   We prove that in general, choosing $\mathbb E$ carefully, we can in fact obtain density in $\mathbb E$.

\begin{lem}\label{dense}
For any $\mathbb D\subseteq G^n$, there is $\mathscr L$-definable $\mathbb B$ such that $\mathbb D$ is a dense subset of $\mathbb B$. Moreover $\mathbb B$ is defined over the same parameters as $\mathbb D$.
\end{lem}

\begin{proof}
Let $\mathbb D$ be definable over $\vec a$. 
By the hypotheses of Theorem \ref{Main Theorem 3}, we know that there are an $\vec a$-definable $\mathbb E$ and $\mathbb S$ such that $\mathbb E$ is  $\mathscr L$-definable, $\mathbb S$ is a dense subset of $G^n$, and $\mathbb D=\mathbb E\cap\mathbb S$. We proceed by induction on the dimension, $k$, of $\mathbb E$ to find an $\mathbb B\subseteq \mathbb E$, $\mathscr L$-definable over $\vec a$ with $\mathbb D$ a dense subset of $\mathbb B$. There is nothing to prove for $k=0$. 

Now suppose we have proven the claim for $j<k$. We may assume that $\mathbb E$ is a cell: write $\mathbb E$ as $\mathbb E_1\cup\dots\cup\mathbb E_l$, with each $\mathbb E_i$ a cell defined over $\vec a$.  If $\mathbb E_i$ is of dimension less than $k$, then we may apply the inductive hypothesis to $\mathbb E_i\cap\mathbb S$. Thus we may assume $\mathbb E$ is a cell of dimension $k$. 

As $\mathbb E$ is a cell, we may choose a projection $\pi:\mathbb E\to \pi(\mathbb E)\subseteq R^k$ so that $\pi$ is a homeomorphism.  Now choose an $\vec a$-definable $\mathbb E'$ and $\mathbb S'$ such that $\mathbb E'$ is  $\mathscr L$-definable, $\mathbb S'$ is a dense subset of $G^k$, and $\pi(\mathbb D)=\mathbb E'\cap\mathbb S'$.  Again, we may divide $\mathbb E'$ into cells, say $\mathbb E'_1\cup\dots\cup\mathbb E'_m$.  For each $i$, either $\mathbb E'_i$ has dimension $k$, in which case it is open and $\pi(\mathbb D)\cap\mathbb E_i'=\mathbb E'_i\cap\mathbb S'$ is dense in $\mathbb E'_i$ or $\mathbb E'_i$ has dimension less than $k$ and we may apply induction to assume $\pi(\mathbb D)\cap\mathbb E_i'$ is a dense subset of $\mathbb E_i'$. Thus $\pi(\mathbb D)$ is a dense subset of $\mathbb E'$. 

Now let $\mathbb B:=\pi^{-1}(\mathbb E')$.  As $\pi$ is a homeomorphism, $\mathbb D$ is a dense subset of $\mathbb B$ and since $\pi$ is $\mathscr L$-definable, so is $\mathbb B$.  We observe that $\mathbb B$ is definable over $\vec a$.
\end{proof}

The second of the two lemmas presents a condition under which a set definable in $(R,G)$ is actually an interval.

\begin{lem}\label{continuity argument}
Let $\mathbb B\subseteq R^n$ be a cell such that $f_1$ and $f_2$ are continuous on $\mathbb B$, $\mathbb B\cap G^n$ is dense in $\mathbb B$, and $f_1(\vec x)<f_2(\vec x)$. Then $\bigcup_{\vec g\in \mathbb B\cap G^n}(f_1(\vec g),f_2(\vec g))$ is an interval.
\end{lem}

\begin{proof}
Let $a=\inf f_1(\mathbb B)$ and $b=\sup f_2(\mathbb B)$.  Let $d\in(a,b)$; we wish to show that $d\in \bigcup_{\vec g\in \mathbb B\cap G^n}\big(f_1(\vec g),f_2(\vec g)\big)$.  For some $c_1\in \mathbb B$, $f_1(c_1)<d$.  Clearly if $f_2(c_1)>d$, we are done, so we may assume that $f_2(c_1)<d$.  Likewise we may assume that there is some $c_2$ such that $d<f_1(c_2)<f_2(c_2)$.  Note that $(f_1+f_2)(c_1)<2d$ while $(f_1+f_2)(c_2)>2d$.  Thus, by the continuity of $f_1$ and $f_2$, and by the connectedness of $\mathbb B$, there is $c_3$ such that $(f_1+f_2)(c_3)=2d$.  Since $f_1<f_2$, we conclude that $d\in(f_1(c_3),f_2(c_3))$.  By the density of $G^n\cap \mathbb B$ in $\mathbb B$ we may find $\vec g\in \mathbb B\cap G^n$ such that $d\in (f_1(\vec g),f_2(\vec g))$. 
\end{proof}

\begin{cor} \label{finite union}
If $f_1$ and $f_2$ are functions $R^n\to R$ which are definable in $\mathscr L$, then $\bigcup_{\vec g\in G^n}(f_1(\vec g),f_2(\vec g))$ is a finite union of intervals.
\end{cor}

\begin{proof}
Recall that $\mathfrak R$ restricted to $\mathscr L$ is an o-minimal structure.  Given $f_1$ and $f_2$, $\mathscr L$-definable $n$-ary functions, we can decompose $R^n$ as a finite union of disjoint cells, $\mathbb C_i$, where on each $\mathbb C_i$, $f_1, f_2$ are continuous, and either the functions coincide on every point of $\mathbb C_i$ or else one of the functions is strictly larger on every point of $\mathbb C_i$.  By Lemma \ref{dense}, we may shrink each $\mathbb C_i$ until we obtain a cell, $\mathbb B_i$, such that $\mathbb B_i\cap G^n$ is a dense subset of $\mathbb B_i$.  By Lemma \ref{continuity argument}, on each such cell, $\bigcup_{\vec g\in \mathbb B_i\cap G^n}(f_1(\vec g),f_2(\vec g))$ is an interval.
\end{proof}

\begin{prop}\label{partition}
Let $\mathbb D\subseteq R$ be definable in $\mathfrak R$. 
Then there is a finite partition $-\infty=a_0<a_1<\dots<a_m=\infty$ of $R$ such that $\mathbb D$ is either small or cosmall in $(a_{i-1},a_i)$ for $i=1,\dots,m$.  Furthermore, if $\mathbb D$ is definable from $\vec d$, so is the partition $-\infty=a_0<a_1<\dots<a_m=\infty$.
\end{prop}
\begin{proof}
We first assume that $\mathbb D$ is basic. So $\mathbb D=\bigcup_{\vec g\in G^n}\varphi(R,\vec g)$, where $\varphi(x,\vec y)$ is an $\mathscr L$-formula. By the o-minimality of $\mathfrak{R}|_{\mathscr L}$, each $\varphi(x,\vec g)$ defines a finite union of points and intervals, and there is a uniform bound on the number of these points and intervals. By Fact \ref{fact 1}, we may assume without loss of generality that each $\varphi(x,\vec{g})$ defines either a single point or a single interval. 

First let us consider the case where $\varphi(x,\vec{g})$ is a single point.  As there is a definable surjection from $G^n$ onto $\mathbb D$, we see that $\mathbb D$ is small. 

Now we consider the case where each $\varphi(x,\vec{g})$ is an interval.  There are $\mathscr L$-definable $f_1,f_2\colon R^n\to R$ such that $\varphi(R,\vec g)=(f_1(\vec g),f_2(\vec g))$.    By Corollary \ref{finite union},
$$\bigcup_{\vec g \in G^n}(f_1(\vec g),f_2(\vec g))$$
is a finite union of intervals. By o-minimality, the endpoints of these intervals are definable over any parameters from which the finite union of intervals may be defined.   

Thus, we have our result if $\mathbb D=\bigcup_{\vec g\in G^n}\varphi(R,\vec g)$.

Now assume $\mathbb D$ and $\mathbb E$ satisfy the conclusion.  To complete the proof, we must show that $\mathbb D^c$ and $\mathbb D\cup\mathbb E$ also have the desired property. But this is clear.
\end{proof}

\begin{defn}\label{small closure}
We say that $\vec e$ is in the {\it small closure} of $A$ iff $\vec e$ is contained in a small set defined with parameters from $A$. 
We denote the small closure of $A$ by $\scl(A)$. 
\end{defn}

\begin{defn}
We say that a set, $\mathbb{S}\subset R^k$, is $G$-bound iff there is an $\mathscr L$-definable $f:R^n\to R^k$ such that $\mathbb{S}\subseteq f(G^n)$.
\end{defn}

It is clear that $G$-bound implies small.  We proceed to prove the converse.

\begin{lem}
Any basic small set $\mathbb{S}$ is $G$-bound.  Furthermore, assuming that there are two definable elements of $R$, the function $f$ witnessing that $\mathbb S$ is $G$-bound is definable over the same parameters as $\mathbb S$.
\end{lem}
\begin{proof}
Note that if $\mathbb S\subset R^k$ is a basic small set, so is each projection of $\mathbb S$ to $R$; and that cartesian products of $G$-bound sets are $G$-bound.  Thus it suffices to consider small subsets of $R$.

Suppose that $\mathbb{S}$ is defined with parameters $\vec a$.  Let $\mathbb{S}$ be
 $$\bigcup_{\vec{g}\in G^n}\varphi(R,\vec{g},\vec{a})$$ where
$\varphi(x,\vec{y},\vec{z})$ is a parameter-free $\mathscr L$-formula. Since  $\mathfrak{R}|_\mathscr{L}$ is o-minimal, each set
$\varphi(R,\vec{g},\vec{a})$ is a finite collection of points and intervals.
It is easy to see that any set containing an open interval is large, so each $\varphi(R,\vec{g},\vec{a})$ is a finite set. By o-minimality, there is a uniform bound $k$ to the size of $\varphi(R,\vec{g},\vec{a})$ for each $\vec g \in G^n$.

Thus mapping $\vec{g}$ to
$\varphi(R,\vec{g},\vec{a})$ gives us a $k$-valued (and $\vec a$-definable)  
function, $F$, in the language $\mathscr L$ such that $F(G^n)=\mathbb S$. By \ref{many to one}, we may replace this with an
actual function, $f$. (Although if $0$ is the only definable element of $R$, we may have to add an additional parameter in $R$.) 
\end{proof}

\begin{rem}
Note that even when $0$ is the only definable element, $\mathbb S$ is still the image  $G$ under a $k$-valued function which is definable with the same parameters as $\mathbb S$.
\end{rem}

\begin{lem} \label{basic small sets rock}
Let $\varphi(x,\vec d)$ define $\mathbb D$. Then there are a partition $-\infty=a_0<\dots<a_n=\infty$ and basic small sets $\mathbb S_1,\dots,\mathbb S_n$ such that $\mathbb D\cap [a_{i-1},a_i]$ either is contained in $\mathbb S_i$, or contains $\mathbb S_i^c\cap[a_{i-1},a_i]$.  Furthermore, the partition and each $\mathbb S_i$ are definable over $\vec d$. 
\end{lem}

\begin{proof}
Note that $\varphi(x,\vec d)$ is equivalent to a boolean combination of basic formulas.  We proceed by induction, using repeatedly that the intersection of a basic set with an interval is again a basic set.

Suppose that $\varphi(x,\vec d)$ is a basic formula. By Proposition \ref{partition}, there is a $\vec d$-definable partition $-\infty=a_0<\dots<a_n=\infty$  such that $\mathbb D\cap [a_{i-1},a_i]$ either is small or cosmall.  If $\mathbb D\cap [a_{i-1},a_i]$ is small, let $\mathbb S_i:=\mathbb D\cap [a_{i-1},a_i]$.  If $\mathbb D\cap [a_{i-1},a_i]$ is cosmall in $[a_{i-1},a_i]$, then $\mathbb D\cap[a_{i-1},a_i]$ is a finite union of intervals, by Lemma \ref{finite union}.  Thus, since it is small, $[a_{i-1},a_i]\setminus\mathbb D$ is a finite collection of points. Let $\mathbb S_i$ be this finite collection of points. Note that in either case, by Proposition \ref{partition},  $\mathbb S_i$ can be defined over $\vec{d}$. 

Now suppose that $\varphi=\varphi_1\land\varphi_2$. Let $\mathbb E_1:=\varphi_1(R,\vec d)$ and let $\mathbb E_2:=\varphi_2(R,\vec d)$.  By induction, there are a partition $-\infty=b_0<\dots<b_m=\infty$  and basic small sets $\widetilde{\mathbb S}_1,\dots,\widetilde{\mathbb S}_m$ with the desired property with respect to $\mathbb E_1$.   Likewise there are a partition  $-\infty=c_0<\dots<c_n=\infty$ and basic small sets $\widetilde{\mathbb S}_{m+1},\dots,\widetilde{\mathbb S}_{m+n}$ with the desired property with respect to $\mathbb E_2$.  Let $-\infty=a_0<\dots<a_l=\infty$ be the union of these two partitions.  Then $\mathbb D \cap[a_{i-1},a_i]$ is either small or cosmall.  

If $\mathbb D \cap[a_{i-1},a_i]$ is small, then either $\mathbb E_1$ or $\mathbb E_2$ is small in $[a_{i-1},a_i]$.  Without loss of generality, we may assume it is $\mathbb E_1$.  Note that $[a_{i-1},a_i]$ is contained in $[b_{k-1},b_k]$ for some $k$.  Let $\mathbb S_i:=\widetilde{\mathbb S}_k\cap[a_{i-1},a_i]$.  As $\widetilde {\mathbb S}_k$ is $\vec d$-definable and contains $\mathbb E_1\cap [b_{k-1},b_k]$, we see that $\mathbb S_i$ satisfies the desired properties.

If $\mathbb D \cap[a_{i-1},a_i]$ is cosmall, then both $\mathbb E_1$ and $\mathbb E_2$ are cosmall in $[a_{i-1},a_i]$.  There are $j$, $k$, such that $[a_{i-1}, a_i]\subseteq[b_{j-1},b_j]$ and $[a_{i-1}, a_i]\subseteq[c_{k-1},c_k]$.  Thus,  $\mathbb E_1\cap[a_{i-1},a_i]$ contains $\widetilde{\mathbb S}_j^c\cap[a_{i-1},a_i]$, and $\mathbb E_2\cap[a_{i-1},a_i]$ contains $\widetilde{\mathbb S}_{m+k}^c\cap[a_{i-1},a_i]$.  Thus, $\mathbb D$ contains $(\widetilde{\mathbb S}_j\cup \widetilde{\mathbb S}_{m+k})^c\cap[a_{i-1},a_i]$.  We let $\mathbb S_i:=(\widetilde{\mathbb S}_j\cup \widetilde{\mathbb S}_{m+k})^c\cap[a_{i-1},a_i]$

Now suppose that $\varphi=\neg\varphi_0$. Let $\mathbb E$ be defined by $\varphi_0$. By induction there is a partition $-\infty=a_0<\dots<a_n=\infty$ and basic small sets $\mathbb S_1,\dots,\mathbb S_n$ such that $\mathbb E\cap [a_{i-1},a_i]$ either is contained in $\mathbb S_i$, or contains $\mathbb S_i^c\cap[a_{i-1},a_i]$, and the $\mathbb S_i$ are defined from $\vec d$. But this partition and these small sets work for $\mathbb D$ as well.

\end{proof}

From the previous two lemmas (as well as Lemma \ref{small in group}), we obtain the following two corollaries: 

\begin{cor}\label{small implies G-bound}
If $\mathbb S$ is a small set, then it is contained in a basic small set and, hence,  $\mathbb S$ is $G$-bound.
\end{cor}
\begin{proof}
 Let $\mathbb S\subset R^k$.  Let $\pi_i$ be the projection onto the $i$th coordinate.  Let $\mathbb S_i:=\pi_i(\mathbb S)$.  By Lemma \ref{basic small sets rock}, take $\widetilde{\mathbb S}_i$, a basic small set containing $\mathbb S_i$.  Then $\widetilde{\mathbb S}_1\times\dots\times\widetilde{\mathbb S}_k$ is a basic small set containing $\mathbb S$.  As $\mathbb S$ is contained in a $G$-bound set, it is itself $G$-bound.
\end{proof}

\begin{cor}\label{small closure is transitive}
A tuple, $\vec e$, is in the small closure of $A$ if and only if there is an $\mathscr L_A$-definable $k$-valued function, $F(\vec x)$ and some $\vec g\in G^n$ such that $\vec e\in F(\vec g)$.  Thus, if $\vec a\in\scl(\vec b)$ and $\vec b\in\scl(\vec c)$  then $\vec a\in\scl(\vec c)$.
\end{cor}

\begin{proof}
If $\vec e\in\scl(A)$ then there is a small set $\mathbb S_{\vec a}$ defined with parameters $\vec a$ from  $A$ that contains $\vec{e}$. The set $\mathbb S_{\vec a}$ is contained in a basic small set, also defined over $A$, and this basic small set is the image of a $k$-valued function on $G^n$.  Conversely, such a set is $G$-bound, and hence small.  Moveover, if $\vec a\in\scl(\vec b)$ and $\vec b\in\scl(\vec c)$ then this is witnessed by $k_1$ and $k_2$-valued functions, $F_1$ and $F_2$ respectively, with $F_1=F_1(\vec x,\vec b)$ and $F_2=F_2(\vec y, \vec c)$.  Thus $F_3:=F_1(\vec x, F_2(\vec y,\vec c))$ witnesses that $\vec a\in\scl(\vec c)$.
\end{proof}

\noindent In addition, we have the following corollary:

\begin{cor}\label{finite union of small sets is small}
A finite union of small sets is again a small set.
\end{cor}

\begin{proof}
Finite unions of $G$-bound sets are again $G$-bound, by Lemma 2.2 of \cite{Mann property}.  
\end{proof}

While we rely on  \cite{Mann property} for the above proof, we note that the corollary also follows as a special case of Proposition \ref{small union of small sets is small} below.

\begin{rem}
Since $\scl$ is transitive, and $\scl(\emptyset)$ is infinite, (and in particular, contains at least one non-zero element) we may add an element of $\scl(\emptyset)$ to the language without affecting small closure.  Thus we may assume that $R$ contains at least two definable elements, and henceforth, we will assume that we may replace each $k$-valued function with an actual function.
\end{rem}

\begin{rem}
Note that, unlike the algebraic closure of $A$, $\scl(A)$ depends on the
model containing $A$.
\end{rem}

Although the following proposition is not used in the proofs of this article's main theorems, it is interesting to note that a small definable union of small sets is again a small set.

\begin{prop}\label{small union of small sets is small}
If $\mathbb{D}$ is small, and $\mathbb{E}_{\vec d}$ is small for each $\vec{d}\in
\mathbb{D}$, then $\bigcup_{\vec d\in\mathbb{D}} \mathbb{E}_{\vec d}$ is also small.
\end{prop}

\begin{proof}
First note that by Corollary \ref{small implies G-bound}, for each $\vec d \in \mathbb{D}$ there is a basic small set containing $\mathbb{E}_{\vec d}$ .  By compactness, the formula defining the basic small set may be chosen uniformly in $\vec d$.  Thus, we may reduce to the case where $\mathbb D$ and each $\mathbb{E}_{\vec d}$ are basic small. 

Assume that the formula $\theta(\vec x,\vec d)$ defines $ \mathbb{E}_{\vec d}$ for every $\vec d \in \mathbb D$. Then, since $\mathbb{E}_{\vec d}$ is a basic small set, there are $\psi(\vec y)\in \tp(\vec d)$ and $f(\vec x,\vec y)$ such that whenever $\vec d'\models \psi(\vec y)$, we have $f(\vec x,\vec d')\colon G^{k}\onto \mathbb{E}_{\vec d'}$. Note that $k$, $\psi$, and $f$ may depend on $\vec d$. However by compactness, there is a finite covering of $\mathbb D$ with sets defined by $\psi_1(\vec y),\dots,\psi_n(\vec y)$, together with associated $k_1,\dots,k_n$ and $f_1,\dots,f_n$. By taking $k=\max\{k_1,\dots,k_n\}$, we see that there is  a definable function $$f(\vec x,\vec y)\colon G^k\times \mathbb D\rightarrow \bigcup_{\vec d\in \mathbb D}\mathbb{E}_{\vec d}$$ such that for any $\vec d\in \mathbb D$, $f(\vec x,\vec d)\colon G^k\onto \mathbb{E}_{\vec d}$.

Now suppose that $g:G^n\onto \mathbb D$ witnesses that $\mathbb D$ is small. Then let $h:G^{k+n}\onto
\bigcup_{\vec d\in\mathbb{D}} \mathbb{E}_{\vec d}$ be defined as follows:
$$h(\vec a_1,\vec a_2):=f(\vec a_1,g(\vec a_2)).$$
So $\bigcup_{\vec d\in\mathbb{D}} \mathbb{E}_{\vec d}$ is $G$-bound, and hence small.
\end{proof}

\begin{defn}
For a set $C$, a function from 
$\mathscr{P}(C)$ to $\mathscr{P}(C)$ is a \emph{closure operator} iff for any $A,B\subseteq C$

(1) $A\subseteq$cl$(A)$,

(2) $A\subseteq B$ implies cl$(A)\subseteq$cl$(B)$,

(3) cl(cl($A$))=cl$(A)$.

\noindent
Furthermore, we say that a closure operator is \emph{finitary} when (2) is strengthened to

(2$'$) $b\in$cl$(A)$
iff $b\in$cl$(A_0)$ for some finite $A_0\subseteq A$.

\noindent
If the closure also satisfies the Steinitz exchange property, then we say
that the closure operator gives rise to a {\it pregeometry}.
\end{defn}

It is clear that the small closure satisfies (1), is finitary, and, by Corollary \ref{small closure is transitive},  satisfies (3). Thus we have proven:
\begin{prop}
The small closure, $\scl$ is a finitary closure operator on subsets of $R$.
\end{prop}

\section{Super-rosiness of $(R,G)$}

In this section we prove Theorem \ref{Main Theorem 3}.  To do this, we will need to use the following propositions from \cite{EO}.  Throughout this section, we assume that $(R,G)$ is $\kappa$-saturated, for $\kappa>2^{\abs{\mathscr L_\mathcal G}}$

\begin{prop}\label{thorn rank and definable maps}
If $\mathbb D$ has \th-rank $\alpha$ and $f:\mathbb D \onto \mathbb E$, then $\mathbb E$ has \th-rank less than or equal to $\alpha$.  Furthermore, if the fibers of $f$ are finite, we have equality.
\end{prop}


\begin{prop} \label{additivity of thorn rank}
If $\mathbb D$ has \th-rank $\alpha$ and $\mathbb E$ has \th-rank less than $\alpha$, then \th-rank$(\mathbb D\setminus\mathbb E)$ is $\alpha$.
\end{prop}

\begin{prop} \label{thorn rank of products}
If $\mathbb D$ has \th-rank $\alpha$ then $\mathbb{D}^n$ has \th-$rank$ at least $\alpha n$, and equality holds if $\alpha=1$.
\end{prop}

Now we begin to analyze \th-dividing in $(R,G)$.  In what follows, $\mathscr L^{\eq}_{\mathcal G}$ refers to the language of $(R,G)^{\eq}$.

\begin{lem}\label{intervals do not thorn divide} Let $\varphi(x,\vec b_0)$ be a formula in $\mathscr L^{\eq}_{\mathcal G}$ with $x$ a variable in the real sort.  If $\varphi(R,\vec b_0)$ is an infinite set definable in $\mathscr L$, then $\varphi(x, \vec b_0)$ does not \th-divide over the empty set. 
\end{lem}

\begin{proof}
It may be worth pointing out that merely because the set $\varphi(R, \vec b_0)$ is definable in $\mathscr L$, we may not assume that $\varphi$ is an $\mathscr L$-formula.  For instance, $\vec b_0$ may come from a sort that does not even exist in $(\mathfrak R|_\mathscr L)^{\eq}$.

Assume, for a contradiction, that $\varphi(x,\vec b_0)$ does \th-divide over the empty set.  That is, $\tp(\vec b_0)$ is non-algebraic, and there is some $\theta(\vec y,\vec c)$ and some $k\in \mathbb N$ such that whenever $\vec b_1,\dots,\vec b_k$ are distinct elements of $\theta(R^{\eq}_{\vec y},\vec c)$, we have that $\varphi(x,\vec b_1)\wedge\dots\wedge\varphi(x,\vec b_k)$ is inconsistent.  Since $\varphi$ defines an infinite $\mathscr L$-definable set, by the o-minimality of $\mathfrak R|_\mathscr L$, it defines a finite collection of points and open intervals.  

First note that we may assume that for each $\vec b\models\theta(\vec y,\vec c)$, it is the case that $\varphi(x,\vec b)$ defines a single interval, modifying $\varphi$ and $\theta$ if necessary.  (It is possible that for some $\vec b\models \theta(\vec y,\vec c)$, $\varphi(x,\vec b)$ defines a finite collection of points.  First we modify $\theta$ to rule out this possibility.  Then we replace $\varphi(x,\vec b)$ with a formula defining the least of the intervals in the finite collection of points and intervals composing $\varphi(R,\vec b)$.)

Now we wish to reduce to the case where $k=2$.  We may assume that $\varphi(x,\vec y)$ does not $(k-1)$-\th-divide. Replace $\varphi(x,\vec y)$ with $$\widetilde{\varphi}(x,\vec y_1,\dots,\vec y_{k-1}):=\bigwedge_{i<k}\varphi(x,\vec y_i)$$ and replace $\theta$ with $$\widetilde{\theta}(\vec y_1,\dots,\vec y_{k-1}):=\theta(y_1)\land\dots\land\theta(y_{k-1})\land\bigwedge_{i<j<k}\vec y_i<\vec y_k.$$  Now $\widetilde{\varphi}$ clearly $2$-\th-divides.

Now we would like to find a contradiction by considering the union of the sets defined by $\varphi(x,\vec b)$ for $\vec b\models\theta$, intersecting with $G$, and noting that it violates (3) of our assumptions on $\mathfrak R$ from Theorem \ref{Main Theorem 3}.  First note that since $G$ is a dense subset of $(a,\infty)$, we can assume that $\varphi(\mathfrak R,b)$ is contained in the closure of $G$ for each $b\models\theta$ (possibly after reflecting the whole family over $a$ and modifying $\theta$.  However, there is still no immediate contradiction since $\bigcup_{\vec b\models\theta}\varphi(R,\vec b)\cap G$ might still be a finite union of intervals in $G$. We can modify $\varphi(x,\vec b)$ once again to define the interval with half the length but the same center as $\varphi(x,\vec b)$.  Now, the union of these intersect $G$ cannot be written as a finite union of intervals intersect a dense subset of $G$.  

\end{proof}

\smallskip

\noindent Now we have all the tools in place to begin our proof of Theorem \ref{Main Theorem 3}.

\medskip

\noindent\textbf{Theorem \ref{Main Theorem 3}.}
\emph{$\mathfrak R=(R,G)$ is super-rosy of \th-rank less than or equal to $\omega$ and \th-rank of $G$ is 1,  Moreover, if $\mathfrak R$ includes a field structure, \th-rank of $\mathfrak R$  equals $\omega$. }

\smallskip

\begin{proof}
First we wish to show that the \th-rank of $G$ is $1$.  For a contradiction, suppose that some formula $\varphi(x,\vec b)$ which defines an infinite subset of $G$ \th-divides over the empty set.  Say that $k$, $\theta(\vec y,\vec c)$ are such that $\bigwedge_{i\leq k}\varphi(x,\vec b_i)$ is inconsistent for any $k$ distinct elements $\vec b_1,\dots,\vec b_k$ satisfying $\theta(\vec y,\vec c)$.  

Then, by $(3)$ of the hypotheses of Theorem \ref{Main Theorem 3}, $\varphi(R,\vec b)$ is a finite union of sets, each of which is either a point or an interval intersect an $\emptyset$-definable dense subset of $G$.  Without loss of generality, we may assume that for each $\vec b'\models \theta(\vec y, \vec c)$, it is the case that $\varphi(x,\vec b')$ defines a single interval, $\psi_1(R,\vec b')$, intersect an $\emptyset$-definable dense subset of $G$.  Which $\emptyset$-definable set may depend on the type of $\vec b'$, but one such set, $\psi_2(R)$, must occur for infinitely many $\vec b'$.  Modifying $\theta$ if necessary, we may assume that for all $\vec b'\models \theta(\vec y, \vec c)$, we have that $\varphi(x,\vec b')$ defines the same set as $\psi_1(x,\vec b')\wedge\psi_2(x)$.

Thus we have that $\{\psi_1(x,\vec b')\wedge\psi_2(x):\vec b'\models \theta(\vec y,\vec c)\}$ is $k$-inconsistent.  But by Lemma \ref{intervals do not thorn divide}, $\psi_1(x,\vec b')$ does not \th-divide, and so we may find an infinite $B=\{b_i:b_i\models\tp(\vec b/\vec c), i<\alpha\}$  such that $\bigcap_{\vec b_i\in B} \psi_1(R,\vec b_i)$ is nonempty and, hence, contains an open interval $(d_1, d_2)$.  But since $\psi_2(x)$ is a dense subset of $G$, $$\bigcap_{\vec b_i\in B} \varphi(R,\vec b_i)\supseteq (d_1, d_2)\cap\psi_2(R)\neq\emptyset,$$ which is a contradiction.

Second, we wish to show that the \th-rank of $x=x$ is no larger than $\omega$.  Suppose that $\varphi(x,\vec b)$ $k$-\th-divides over the empty set, where, again, $\vec b$ may come from any sort in $\mathfrak{R}^{\textrm{eq}}$.   We observe that it suffices to show that $\mathbb{D}_{\vec b}:=\varphi(R,\vec b)$ must be a small set, since any small set is $G$-bound, and thus we may apply  Proposition \ref{thorn rank and definable maps} and Proposition \ref{thorn rank of products} to conclude that any $G$-bound set has finite \th-rank.  Then we will have shown that any formula, $\varphi(x,\vec b)$, which \th-divides has finite \th-rank, and, thus, $\tr{x=x}\leq\omega$.

Now assume for a contradiction that $\varphi(x,\vec b)$ is not a small set.  By \ref{basic small sets rock} there is some open interval $\mathbb{I}_{\vec b}$ such that $\mathbb{D}_{\vec b}$ is cosmall in $\mathbb{I}_{\vec b}$, that is, $\mathbb D_{\vec b}\cap\mathbb I_{\vec b}=\mathbb I_{\vec b}\setminus \mathbb S_{\vec b}$ where $\mathbb S_{\vec b}$ is a small set.  Suppose that $\theta(\vec y,\vec c)$ is such that for any $\vec{b}_1,\dots,\vec{b}_{k}$, each realizing  $\theta(\vec y,\vec c)$, one has 
$$\mathbb{D}_{\vec{b}_1}\cap\dots\cap\mathbb{D}_{\vec{b}_{k}}=\emptyset.$$  

Thus we have
$$\emptyset= \bigcap_{1\leq i\leq k} (\mathbb D_{\vec b_i}\cap\mathbb I_{\vec b_i})=\bigcap_{1\leq i\leq k} \mathbb I_{\vec b_i} \setminus \bigcup_{1\leq i\leq k} \mathbb S_{\vec b_i}$$

Then it is not hard to see that  
$$\mathbb J:=\mathbb{I}_{\vec{b}_1}\cap\dots\cap\mathbb{I}_{\vec{b}_{k}}=\emptyset.$$

For if this were not the case, $\mathbb J$ would be an open interval contained in the small set $\mathbb S_{\vec b_1}\cup\dots\cup\mathbb S_{\vec b_k}$, which is impossible, by Corollary \ref{finite union of small sets is small}.

Thus, if $\psi(x, \vec b)$ defines $\mathbb I_{\vec b}$, we see that $\psi(x,\vec b)$ also \th-divides.  But since intervals are $\mathscr L$-definable, this contradicts the previous lemma.  Thus we conclude that $\tr{x=x}$ is no greater than $\omega$.

It remains to show that if $R$ has a field structure, then $\tr{x=x}$ is precisely $\omega$.  Note that as $G$ is small, $R$ is an infinite dimensional $\dcl(G)$-vector space.  Choose $(c_i)_{i\in\mathbb{N}}$ independent vectors.  Considering $$c_1 G+\dots+c_{n-1} G+c_n g,$$ and noting that one gets $2$-inconsistency as one varies $g$ though $G$, it is clear that $$\mathbb{V}^n_{\vec c}:=c_1 G+\dots+c_{n-1} G+c_n G$$ has \th-rank $n$.  As each $\mathbb{V}^n_{\vec c}$ is a subset of $R$, $\tr{R}\geq\omega$.
\end{proof}

Note that we have not only shown that $\mathfrak R$ is super-rosy, but the following:

\begin{cor}\label{thorn dividing implies small}
Any formula $\varphi(x,\vec b)$ that \th-divides defines a small subset of $R$.
\end{cor}

This will allow us to show that, in certain cases, small closure gives rise to a pregeometry in Section \ref{thorn-U-rank}.

Finally, we should point out the following two corollaries:
\begin{cor}
    Dense pairs of o-minimal structures (with at least a group structure) are superrosy.  If the o-minimal structure is an expansion of a real closed field, the \th-rank of the pair is $\omega$.
\end{cor}
\begin{proof}
    See \cite{Dense Pairs} to see a proof that dense pairs satisfy the hypotheses of Theorem \ref{Main Theorem 3}.
\end{proof}

For the next corollary, we need a defintion and a fact from \cite{o-minimal open core}:
\begin{defn}
    An expansion of $(\mathbb R,<)$ is said to have \emph{o-minimal open core} if the reduct generated by the definable open sets is o-minimal. 
\end{defn}
\begin{fact}
    An expansion of $(\mathbb R, +, \cdot)$ has o-minimal open core if and only if each definable open subset of $\mathbb R$ has finitely many connected components.
\end{fact}
\begin{cor}\label{o-minimal open core}
    An expansion of $(\mathbb R, +, \cdot)$ which satisfies the hypotheses of Theorem \ref{Main Theorem 3} has o-minimal open core.
\end{cor}
\begin{proof}
    For a contradiction, let $\mathbb D$ be definable, open, and with infinitely many connected components.  We may assume that $\mathbb D\subset(a,\infty)$.  We note that that given $d\in\mathbb D$, the connected component of $\mathbb D$ containing $d$ is definable, say by $\varphi(x,d)$.  Being in the same connected component is a definable equivalence relation, call it $E$.  Thus the connected component of $d$ may just as easily be defined by $\widetilde\varphi(x,d/E)$. As $d/E$ varies through the sort $\mathbb D/E$, $\widetilde\varphi(x,d/E)$ \th-divides.  But $\widetilde\varphi(x,d/E)$ is an interval, and hence $\mathscr L$-definable.  This contradicts Lemma \ref{intervals do not thorn divide}.
\end{proof}

\section{Imaginaries}\label{imaginaries}

Pillay, building on ideas of Lascar, showed that a strongly minimal theory where the algebraic closure of the empty set is infinite eliminates imaginaries down to finite sets (see e.g. \cite{Model Theory of Fields}).  What follows is the same argument, with small replacing finite, and it shows that $\mathfrak{R}$ eliminates imaginaries down to small sets.

In this section, we assume that $(R,G)$ satisfies all the hypotheses of Theorem \ref{Main Theorem 4}.  That is, we add to the assumptions of the last section, the assumption that given any set $A$, and $I$ any interval defined over $A$, that $\scl(A)\cap I$ is not contained in any small set.

\begin{prop}
Let $\varphi(\vec x, \vec y)$ define an equivalence relation, $E$, and let $e$ be an element of the sort $R_{\vec x}/E$. Then there is an element, $\vec d$, of $R_{\vec x}$ such that $e=\vec d/E$ and $\vec d\in\scl(e)$.
\end{prop}

\begin{proof}
Let $\pi:R^n\to R^n/E$ be the quotient map, and consider $\mathbb D_1$ defined by $$\exists x_2,\dots,x_n \pi(x_1,x_2,\dots,x_n)=e.$$  In the case that $\mathbb D_1$ is small,  any element of $\mathbb D_1$ is in $\scl(e)$; let $d_1$ be any such element.  Otherwise, there is some interval such that $\mathbb D_1$ is cosmall in that interval.  By our assumption on the small closure, it is not possible that $\scl(e)$ is contained in $\mathbb{D}_1^c$.  Let $d_1$ be some element of $\scl(e)\cap \mathbb D_1$.
  
Proceed inductively and define $\mathbb D_i$ as $$\exists x_{i+1},\dots,x_n \pi(d_1,\dots,d_{i-1},x_i,x_{i+1},\dots,x_n)=e$$ and consider the cases of $\mathbb D_i$ small, or not, as above, to get $\vec d:=(d_1,\dots,d_n)$.  Then $d_i\in\scl(\vec e,d_1,\dots,d_{i-1})$.  By choice of $d_1,\dots,d_{i-1}$, together with the fact that $\scl:\mathscr P(R) \to \mathscr P(R)$ is a closure operator, this implies that $d_i\in\scl(e)$.

\end{proof}

Now we may prove our elimination of imaginaries result:

\noindent\textbf{Theorem \ref{Main Theorem 4}.}
\emph{Enlarge $\mathfrak R$ to $\widetilde{\mathfrak R}$ by adding sufficiently many sorts of $\mathfrak{R}^{\eq}$ so that $\widetilde{\mathfrak{R}}$ has a code for every basic small subset of $R^k$. Then $\widetilde{\mathfrak{R}}$ eliminates imaginaries.}

\smallskip

\begin{proof}
Take $e\in R^{\eq}$.  We want to find $c\in\widetilde{R}$ such that $c$ is interdefinable with $e$. Take $\vec d$ such that $\pi(\vec d)=e$ and $\vec d\in \scl(e)$.   Thus $\vec d$ is in a basic small set, $\mathbb D$, defined over $e$; let $c$ be the code for $\mathbb D\cap\pi^{-1}(e)$. Clearly, $c$ is defined over $e$.  But $e$ is defined over any element of $\mathbb D\cap\pi^{-1}(e)$, and thus over $c$ as well.

\end{proof}

\section{Groups with the Mann Property} \label{bounded case}

We start by defining the Mann property for multiplicative subgroups of fields. Let $K$ be a field, and $G$ a subgroup of $K^\times$. For $a_1,\dots,a_n\in K$, a solution $(g_1,\dots,g_n)$ of $a_1x_1+\cdots+a_nx_n=1$ in $G$ is said to be {\it nondegenerate} if $\sum_{i\in I}a_ig_i\neq 0$ for every non-empty subset $I$ of $\{1,\dots,n\}$. We say $G$ has the {\it Mann property} if for every $a_1,\dots,a_n$ from  $K$, the equation $a_1x_1+\cdots+a_nx_n=1$ has finitely many nondegenerate solutions in $G$.

\medskip
Prior to this section, we have assumed that $(R,G)$ was as in Theorem \ref{Main Theorem 3}.  In this section we instead prove that $(R,G)$ as in Theorem \ref{Main Theorem 1} satisfy the hypotheses of Theorem \ref{Main Theorem 3}.  That is, we assume that $R$ is a real closed field and $G$ is a dense subgroup of $R^{>0}$ with the Mann property and such that for each $p$, the $p$th powers in $G$ have finite index in $G$. 

As noted in the introduction, most of the results about groups with the Mann property that we need are found in \cite{Mann property}.  For instance, we have the following:

\begin{fact}
 By of Lemma 6.1 of \cite{Mann property}, if $(R,G)$ satisfies the conditions of Theorem \ref{Main Theorem 1}, then $G$ is small.
\end{fact}

\begin{fact}
By Theorem 7.5 of \cite{Mann property}, if $(R,G)$ satisfies the conditions of Theorem \ref{Main Theorem 1}, then any definable subset of $\mathfrak{R}$ is a boolean combination of basic sets.
\end{fact}

\noindent However, we will need to strengthen the quantifier elimination results obtained there. 

\medskip
In the rest of this section $q$ is of the form $p^m$, where $p$ is a prime number and $m\in\mathbb N$.

For each $q$ and $\vec k=(k_1,\dots,k_n)\in\Z^{n}$ let $D_{q,\vec k}(\vec x)$ be the formula
$$\mathcal G(x_1)\land\dots\land \mathcal G(x_n)\land\exists y (\mathcal G(y)\land x_1^{k_1}\cdots x_n^{k_n}=y^q).$$
Note that $D_{q,(0,\dots,0)}(R_{\vec x})$ is all of $G^n$, and for any $g\in G$, there is $\vec h\in G^n$ such that $D_{q,1\vec k}(g,R_{\vec x})$ equals $\vec h D_{q,\vec k}(R_{\vec x})$.

We will write $G^{[n]}$ to denote the elements of $G$ that have $n$th roots in $G$.

\medskip

\begin{prop}\label{induced structure}
Let $\mathbb D\subseteq G^n$ be definable in $(R,G)$, then $\mathbb D$ is a boolean combination of sets of the form $\mathbb F\cap\vec g D_{q,\vec k}(R_{\vec x})$, where $\mathbb F$ is a semialgebraic set, $\vec g\in G^n$, $q$ is as above, and $\vec k\in\Z^n$.
\end{prop}
\medskip
Before proving this proposition, we recall some results from \cite{Mann property} that are used in the proof of it.

Let $(R_1,G_1)$ and $(R_2,G_2)$ be two $|R|^+$-saturated elementary extensions of $(R,G)$. Then in the proof of Theorem 7.1 of \cite{Mann property}, the authors construct a back and forth system $\mathcal I$, between $(R_1,G_1)$ and $(R_2,G_2)$, consisting of isomorphisms $\iota:(R_1',G_1')\to(R_2',G_2')$ where
$R_i'$ is a real closed ordered subfield of $R_i$ of cardinality $<|R|$,
$G_i'\subseteq R_i'^{>0}$ is a pure subgroup of $G_i$ containing $G$, and 
$R_i'$ and $\Q(G_i)$ are algebraically free over $\Q(G_i')$ for $i=1,2$. 
\smallskip

We also need the following lemma from \cite{Mann property}.

\begin{lem}\label{independence}
Let $R$ be a real closed field with a subfield $E$ and let $H\subseteq R^{>0}$ be a subgroup satisfying the Mann property. Suppose that $H'$ is a subgroup of $H$ such that for all
$a_1,\dots,a_n\in E^{\times}$ the equation 
$a_1x_1 + \dots + a_nx_n=1$ has the same nondegenerate solutions
in $H'$ as in $H$. Then for any $h\in H$, if $h$ is algebraic over 
$E(H')$ of degree $d$, then $h^d\in H'$.
\end{lem}

\medskip\noindent
Now we prove Proposition \ref{induced structure}.

\begin{proof}
By standard model theoretic arguments (see for instance 8.4.1 of \cite{Hodges}), it is enough to prove the following:

\smallskip
\noindent{\bf Claim.} Let $(R_1,G_1)$ and $(R_2,G_2)$ be two $|R|^+$-saturated elementary extensions of $(R,G)$. Take $\vec g_1\in G_1^n$ and $\vec g_2\in G_2^n$ such that for any formula $\varphi(\vec x)$ in the language of ordered rings with parameters in $R$, for any $g\in G$, and for any $q$, $\vec k$ as above, we have 
$$(R_1,G_1)\models\varphi(\vec g_1)\land D_{q,1\vec k}(g,\vec g_1)\textrm{ iff }(R_2,G_2)\models\varphi(\vec g_2)\land D_{q,1\vec k}(g,\vec g_2).$$ 
Then $(R_1,G_1,\vec g_1)\equiv_{R}(R_2,G_2,\vec g_2)$. 
\smallskip

\noindent{\it Proof of the claim.} By the remarks made before the proof, there is a back and forth system $\mathcal I$ between $(R_1,G_1)$ and $(R_2,G_2)$. It suffices to prove that there is an element $\iota$ of $\mathcal I$ taking $\vec g_1$ to $\vec g_2$. 

Since $\vec g_1$ and $\vec g_2$ satisfy the same ordered field type over $R$, there is a ordered field isomorphism $\iota:R_1'\to R_2'$, mapping $\vec g_1$ to $\vec g_2$ equal to the identity on $R$, where $R_i'$ is the real closure of $R(\vec g_i)$ for $i=1,2$.  

Consider $G_i':=R_i'\cap G_i$. We wish to show that $G_i'=G\<\vec g_i\>:=\{(g\vec g_i^{\,\vec k})^{1/m}: g\in G,\vec k\in\Z^n,m\in\N, g\vec g_i^{\,\vec k}\in G_i^{[m]}\}$. It is clear that $G_i'\supseteq G\<\vec g_i\>$. 

We use Lemma \ref{independence} to show $G_i'\subseteq G\<\vec g_i\>$. To do this we need to check that for all $a_1,\dots a_n\in R$, if $a_1 x_1+\dots+a_k x_n=1$ has a nondegenerate solution in $G_i$, then this solution lies in $G\<\vec g_i\>$.  But since $(R,G)\preceq(R_i,G_i)$, such a solution lies even in $G$.  Now applying Lemma \ref{independence}, we see that if $g\in G_i$ is algebraic of degree $d$ over $R(G\<\vec g_i\>)$, then $g^d$ is in $G\<\vec g_i\>$ and thus $g$ itself is in $G\<\vec g_i\>$.

Now we wish to show that $\iota(G_1')=G_2'$. An element of $G_1'$ is of the form $(g\vec g_1^{\,\vec k})^{1/m}$ for some $g\in G,\vec k\in\Z^n,m\in\N$. Note $\iota((g\vec g_1^{\,\vec k})^{1/m})=(g\vec g_2^{\,\vec k})^{1/m}$, and by our assumption on $\vec g_i$, $(g\vec g_1^{\,\vec k})$ is in $G_1^{[m]}$ if and only if $(g\vec g_2^{\,\vec k})$ is in $G_2^{[m]}$.  Thus $\iota$ is an isomorphism between $(R_1',G_1')$ and $(R_2',G_2')$.

It remains to show that $R_i'$ and $\Q(G_i)$ are algebraically free over $\Q(G_i')$ and $G_i'$ is a pure subgroup of $G_i$. The first follows from the assumption that $(R_i,G_i)$ is an elementary extension of $(R,G)$, and $G_i'$ is a pure subgroup of $G_i$, since it equals $G\<\vec g_i\>$.
\end{proof}

\begin{rem}\label{union}
Note that the proof of Proposition \ref{induced structure} does not require that the subgroup of $p$th powers has finite index. With this assumption, we see that in addition, the subgroup of $q$th powers is of finite index in $G$ and therefore $D_{q,\vec k}(R_{\vec x})$ is of finite index in $G^n$. So $G^n\setminus D_{q,\vec k}(R_{\vec x})$ is a finite union of cosets of $D_{q,\vec k}(R_{\vec x})$. 
\end{rem}

\medskip\noindent
We also have the following lemma.

\begin{lem}\label{D_{q,n} is dense}
For any $q$, and $\vec k\in\mathbb Z^n$, $D_{q,\vec k}(R_{\vec x})$ is dense in $G^n$.
\end{lem}

\begin{proof}
We show that for any $q$, and $\vec k\in\mathbb Z^n$, $D_{q,\vec k}(R_{\vec x})\supseteq (G^{[q]})^n$, which is enough to prove the lemma, as $(G^{[q]})^n$ is dense in $G^n$. So let $(g_1^q,\dots,g_n^q)\in(G^{[q]})^n$. Then $$(g_1^q)^{k_1}\cdots(g_n^q)^{k_n}=(g_1^{k_1})^q\cdots(g_n^{k_n})^q=(g_1^{k_1}\cdots g_n^{k_n})^q\in G^{[q]}.$$
Thus $(g_1^q,\dots,g_n^q)\in D_{q,\vec k}(R_{\vec x})$.
\end{proof}

\begin{cor}\label{Lou's Corollary}
    Each $D_{q,\vec k}(R_{\vec x})$ is a finite union of cosets of $(G^{[q]})^n$.  Moreover, for any $\mathbb D\subset G^n$ there is $d\in \mathbb N$ such that $\mathbb D$ is a finite union of sets of the form $\mathbb F\cap \vec g(G^{[d]})^n$ where $\mathbb F$ is semialgebraic.\footnote{The authors thank Lou van den Dries for pointing out this Corollary.}
\end{cor}

\begin{proof}
    By the proof of Lemma \ref{D_{q,n} is dense}, we have that $(G^{[q]})^n$ is a subgroup of $D_{q,\vec k}(R_{\vec x})$. Since $(G^{[q]})^n$ is finite index in $G^n$, it is also finite index in $D_{q,\vec k}(R_{\vec x})$.
    
    Next note that if $d$ is the least common multiple of $d_1,d_2$, then $G^{[d_1]}\cap G^{[d_2]}=G^{[d]}$.  Thus, given any finite number of cosets of $(G^{[d_i]})^n$ for various $d_i$, one may replace them by a finite number of cosets of $(G^{[d]})^n$, where $d$ is the least common multiple of the $d_i$.  Using this observation, the reader may easily check that for each $\mathbb D\subset G^n$ there is $d\in \mathbb N$ such that $\mathbb D$ is a finite union of sets of the form $\mathbb F\cap \vec g(G^{[d]})^n$ where $\mathbb F$ is semialgebraic.
\end{proof}

Now we are in a position to prove the first of our main results.

\smallskip

\noindent\textbf{Theorem \ref{Main Theorem 1}.}
\emph{$\mathfrak R=(R,G)$ is super-rosy of \th-rank equal to $\omega$ and \th-rank of $G$ is 1. }

\smallskip

\begin{proof}
Since super-rosiness and \th-rank are properties of the theory, we may assume that $(R,G)$ is sufficiently saturated.  Conditions (1) and (2) of Theorem \ref{Main Theorem 3} are clear; we will show (3) for $(R,G)$ in a language expanded by naming each element of some model.   Consider $\mathbb D\subseteq G^n$. First, we wish to show that $\mathbb D=\mathbb E\cap\mathbb S$, where $\mathbb E$ is semialgebraic and $\mathbb S$ is a dense subset of $G^n$. For the purposes of this proof, we refer to such sets as {\em nice}.

We have established, in the previous corollary, that $\mathbb D=\bigcup_{i=1}^m \mathbb E_i\cap \mathbb S_i$, where each $\mathbb E_i$ is semialgebraic, and each $\mathbb S_i$ is of the form $\vec g(G^{[d]})^n$, and, in particular, each $\mathbb S_i$ is dense in $G^n$.  Thus $\mathbb D$ is a finite union of nice sets.  We wish to show that a finite union of nice sets is nice.  Consider $(\mathbb E_1\cap\mathbb S_1) \cup (\mathbb E_2\cap \mathbb S_2)$.  Let $\widetilde{\mathbb E}_1:=\mathbb E_1\setminus \mathbb E_2$ and $\widetilde{\mathbb E}_2:=\mathbb E_2\setminus\mathbb E_1$.  Let $\widetilde{\mathbb S}_1:=\mathbb S_1\setminus \mathbb S_2$ and $\widetilde{\mathbb S}_2:=\mathbb S_2\setminus\mathbb S_1$.  Let $\mathbb E:=(\mathbb E_1 \cup \mathbb E_2)$ and let $\mathbb S=(\mathbb S_1\cup\mathbb S_2)\setminus ((\widetilde {\mathbb E}_2 \cap \widetilde{\mathbb S}_1)\cup (\widetilde{\mathbb E}_1\cap \widetilde{\mathbb S}_2))$.  Note that 
$$(\mathbb E_1\cap\mathbb S_1) \cup (\mathbb E_2\cap\mathbb S_2)=\mathbb E\cap \mathbb S.$$
Thus we want to show that $\mathbb S$ is dense in $G^n$.

Suppose that $\mathbb S$ is not dense in $G^n$.  Then there is an semialgebraic open $\mathbb U\subseteq (R^{>0})^n$ such that $\mathbb S\cap\mathbb U=\emptyset$.  Thus $\mathbb S_1\cap \mathbb U\subseteq\widetilde{\mathbb E}_2\cap \widetilde{\mathbb S}_1$ and $\mathbb S_2\cap \mathbb U\subseteq\widetilde{\mathbb E}_1\cap \widetilde{\mathbb S}_2$. Since $\mathbb S_2$ is dense,  the closure of $\mathbb S_2\cap \mathbb U$  equals the closure of $\mathbb U$, and is contained in the closure of $\widetilde{\mathbb E}_1$.  Thus, $\widetilde{\mathbb E}_1$ must contain all of $\mathbb U$ except for a semialgebraic set, $\mathbb D_1$, of dimension less than $n$. Likewise there is $\mathbb D_2$ such that $\mathbb U\setminus \mathbb D_2 \subseteq \widetilde{\mathbb E}_2$.  Thus $\mathbb U\setminus (\mathbb D_1 \cup \mathbb D_2)$ is contained in $\widetilde{\mathbb E}_1\cap \widetilde{\mathbb E}_2=\emptyset$, a contradiction.

Finally we note that by Corollary \ref{Lou's Corollary}, if $\mathbb D\subseteq G^n$, then $\mathbb D=\bigcup_{i<k}(\mathbb E_i\cap\mathbb S_i)$ with each $\mathbb E_i$ a semialgebraic set and each $\mathbb S_i$ a coset of $(G^{[d]})^n$.  Since $(G^{[d]})^n$ is a subgroup of finite index, any model has representatives of each coset, and thus, after naming the elements any model, each $\mathbb S_i$ becomes $\emptyset$-definable, and we may apply Theorem \ref{Main Theorem 3} to get that $(R,G)$ in this expanded language is super-rosy of \th-rank equal to $\omega$ and \th-rank of $G$ is 1.  Since \th-rank is invariant under expansions of the language by constants, we are done.
\end{proof}

In \cite{o-minimal open core}, the question is raised whether $(\mathbb R, 2^\mathbb Z 3^\mathbb Z)$ has o-minimal open core.  We are now in a position to give an affirmative answer to this question.

\begin{cor}
   If $(\mathbb R, G)$ is an expansion of the real field by a predicate for a dense multiplicative subgroup  of $\mathbb R^{>0}$ with the Mann property, then $(\mathbb R,G)$ has o-minimal open core.
\end{cor}
\begin{proof}
    By Corollary \ref{o-minimal open core}.
\end{proof}

To prove the second main result, that adding codes for the small sets definable in $\mathfrak R$ is sufficient to eliminate imaginaries, we must verify our assumptions at the beginning of Section \ref{imaginaries}: that given any set of parameters $A$, and any interval $I$ defined over $A$, the small closure of $A$ intersect $I$ is not contained in any small set.  To do this, we must first perform some \th-rank calculations within $\mathfrak R$.

\begin{defn}
For $n>0$ we define $G^{+n}$ inductively as $$G^{+1}:=G\cup\{0\},$$ $$\text{and }G^{+(n+1)}:=(G\cup\{0\})+G^{+n}.$$

\end{defn}

\begin{prop} \label{thorn rank of the n-fold sum of G}
The \th-rank of $G^{+n}$ is $n$.
\end{prop}

\begin{proof}
Consider the map $f:G^n\to G^{+n}$ given by $f(\vec g)=g_1+\cdots+g_n$. We have \th-rank of $G^{+n}$ is less than or equal to $n$, since $f$ is surjective. 

For the converse, define $G^n_I:=\{\vec g\in G^n: \sum_{i\in I}g_i=0\}$ for any nonempty subset $I$ of $\{1,\dots,n\}$. Note that $G^n_I$ is the image of $G^{n-1}$ under a definable map, thus is of \th-rank at most $n-1$. Now define
$$G^n_{\text{nd}}:=G^n\setminus\bigcup_{\emptyset\neq I\subseteq\{1,\dots,n\}}G^n_I.$$
Note that \th-rank of $G^n_{\text{nd}}$ is $n$, and by the Mann property, the restriction of $f$ to $G_{\text{nd}}^n$ has finite fibers. Therefore, by \ref{thorn rank and definable maps}, \th-rank of $G^{+n}$ is $n$.
\end{proof}

\begin{prop}
Let $A$ be any set, and $I$ any interval defined over $A$. Then $\scl(A)\cap I$
is not contained in any small set.
\end{prop}

\begin{proof}
Note that $\scl(A)$ contains $\scl(\emptyset)$ which in turn contains
$G^{+ n}$. First we show that $$\bigcup_{n>0}
G^{+n}$$ is not contained in any small set. Assume it is contained in a small set $\mathbb S$. Since $\mathbb S$ is $G$-bound, there is a map $f:R^k\to R$ such that $\mathbb S\subseteq f(G^k)$. Therefore by Propositions \ref{thorn rank and definable maps} and \ref{thorn rank of products}, we have \th-rank of $\mathbb S$ is at most $k$, and thus, for each $n$, $G^{+n}$ has \th-rank at most $k$ contradicting Proposition \ref{thorn rank of the n-fold sum of G}.

Let $I=(b,c)$.  Now let $f:R\to(b,c)$ be a definable bijection.  Note that $f(\bigcup_{n>0} G^{+n})$ is contained in $\scl(A)\cap I$.  If $f(\bigcup_{n>0} G^{+n})$ were contained in some small set, say $\mathbb{S}$, then $f^{-1}(\mathbb{S})$ would be a small set containing $\bigcup_{n>0} G^{+n}$, a contradiction.
\end{proof}

\noindent Now we have proven the second of main results:

\smallskip

\noindent\textbf{Theorem \ref{Main Theorem 2}}
\emph{If one enlarges $(R,G)$ by adding sufficiently many sorts of $(R,G)^{\eq}$ so that the resulting structure has a code for every basic small subset of $R^k$, then this structure eliminates imaginaries.}

\section{The structure $R^{>0}/G$}

In this section we assume that $R$ has a field structure.

\begin{prop}\label{humpty-dumpty}
Let $C\subset R$ and let $a$, $b\in R$ be such that $a,b\not \in \scl(C)$. Then for every formula $\varphi(x,\vec c)$ in $\tp(a/C)$ there is $b'\in R$ such that $b'/G=b/G$ and $b'\in \varphi(R,\vec c)$.
\end{prop}

\begin{proof}
We may assume that $C=\dcl(C)$. Let $\varphi(x,\vec c)\in \tp(a/C)$. By Lemma \ref{basic small sets rock} there is a partition $\{c_0,\dots,c_n\}$ of $R$, where $c_i\in C$ for $i\leq n$ such that $\varphi(x,\vec c)$ is small or cosmall when restricted to $(c_i,c_{i+1})$.  Say $a\in (c_i,c_{i+1})$. Since $a\not \in \scl(C)$, $\varphi(R,\vec c)$ is cosmall in $(c_i,c_{i+1})$. Since $b\neq 0$, there is $t\in R$ such that $tb=a$. Furthermore, since multiplication by $b$ is a continuous function, and since $G$ is dense in $R$, we can find $g\in G$ such that $b'=gb\in (c_i,c_{i+1})$. We may choose $g$ \th-independent from $b$ over $C$. Since $b\not \in \scl(C\cup \{g\})$ and multiplication by $g$ is a definable bijection of $R$, we have that $b'\not \in \scl(C\cup\{g\})$ and thus  $\varphi(x,\vec c)\in \tp(b'/C)$.

\end{proof}

\begin{cor}
Let $a$, $b\in R$ be such that $a,b\not \in \scl(A)$. Let $a_G=a/G$, $b_G=b/G$. Then for any set $A$ such that $a_G$ and $b_G$ are \th-independent from $A$, $\tp(a_G/A)=\tp(b_G/A)$.
\end{cor}

\begin{proof}
We may assume that $a$ and $b$ are independent from $A$. By the previous proposition for every formula $\varphi(x,\vec c)$ in $\tp(a/A)$ we can find $b'\in R$ such that $b'/G=b_G$ and $b'\in \varphi(R,\vec c)$. This implies that $\tp(a_G/A)=\tp(b_G/A)$.
\end{proof}

Given any subset $C\subset \mathfrak R$, there is a unique type in $R^{>0}/G$ over $C$ that contains only large sets. Thus the group $R^{>0}/G$ is definably connected (in the sense of having no proper definable subgroups of finite index) and all definable subsets of $R^{>0}/G$ are small or cosmall. 

Assume now that $R$ is uncountable and $G$ is countable. Then the definable small sets are countable.  This raises the following question:

\begin{ques}
Is $R^{>0}/G$ quasi-minimal?
\end{ques}
 
In \cite{quasi-minimal}, Zilber defines a quasi-minimal excellent class, as a class of structures closed under isomorphism, where each definable set is countable or co-countable, and with a closure operator satisfying three assumptions.   When, in addition, the closure operator satisfies the exchange property, he obtains that the class is categorical in every uncountable cardinal.  We have that each definable set is countable or co-countable, and small closure satisfies exchange and can easily be seen to satisfy the first of Zilber's three assumptions.  However, we have been unable to verify that the other two assumptions hold.

Even without the assumption that $G$ is countable, we may ask the following, less ambitious, question:  

\begin{ques}
Is $R^{>0}/G$ superstable?
\end{ques}

There is no obvious order definable within $R^{>0}/G$, and if $R^{>0}/G$ does not have the order property, it must be superstable, as \th-forking agrees with forking in stable theories.

\section{The $\uth$-rank}\label{thorn-U-rank}
Throughout this section, $\mathfrak R$ denotes a structure satisfying the hypotheses of Theorem \ref{Main Theorem 3}.

In \cite{Bu} Buechler used infinite dimensional pairs to study the geometric properties of a strongly minimal sets. He showed the pair has Morley rank one iff the strongly minimal set is trivial, Morley rank two iff the strongly minimal set is locally modular non trivial and $\omega$ otherwise. These results were generalized by Vassiliev in \cite{Va1} to the setting of simple theories using lovely pairs to analyze $SU$ rank one pregeometries. Dense pairs of o-minimal structures were studied by van den Dries in \cite{Dense Pairs}, where he showed they satisfy the hypothesies of Theorem \ref{Main Theorem 3}. In what follows below, we show that the same relationship exists between the pregeometry of a o-minimal  structure, and that of the corresponding dense pair (though, of course, here the information yielded by the dense pair is already known).

Peterzil and Starchenko \cite{PS} showed that locally every o-minimal structure behaves as an expansion of a field, an ordered vector space, or is trivial. In the analysis that follows below, we will deal with two cases: when $\mathfrak R$ includes a field structure and when $\mathfrak R|_\mathscr L$ is an ordered abelian group with no additional structure.

Recall that the $\uth$-rank ``counts'' the number of times the type can \th-fork and that $1$-types in o-minimal structures have $\uth$-rank at most one.

\begin{lem}
Let $g\in G$ and let $C\subset R$. Then $\uth(\tp(g/C))\leq 1$ and equality holds iff $g\not \in \dcl(C)$.
\end{lem}

\begin{proof}
It follows from Theorem \ref{Main Theorem 3}.
\end{proof}

\subsection{Field case}
Now assume that $\mathfrak R|_\mathscr L$ has a definable field structure. Then, as $G$ is small, $R$ is an infinite dimensional $\dcl(G)$-vector space and we fix a countable family $(c_i)_{i \in \omega}$ of linearly independent vectors.

\begin{defn}
Let $g_1,\dots,g_n\in G$ and let $A\subset R$. We say that $\{g_1,\dots,g_n\}$ is an $A$-{\it independent set} if $\uth(\tp(g_1,\dots,g_n/A))=n$.
\end{defn}

\begin{lem} \label{G-bound gives finite U-rank}
Let $g_1,\dots,g_n\in G$ and let $C=\{c_1,\dots,c_n\}$. Then $$\uth(\tp(c_1g_1+\dots+c_ng_n/C))\leq n$$ and equality holds iff $\{g_1,\dots,g_n\}$ is a $C$-independent set.
\end{lem}

\begin{proof}
Clearly $c_1g_1+\dots+c_ng_n\in \dcl(\{g_1,\dots,g_n,c_1,\dots,c_n\})$, so by additivity of the rank and the previous lemma, $$\uth(\tp(c_1g_1+\dots+c_ng_n/C))\leq\uth(\tp(g_1,\dots,g_n/C))\leq n.$$ 
Furthermore since $C=\{c_1,\dots,c_n\}$ is a set of linearly independent vectors, there is only one solution in $G^n$ for the equation $c_1x_1+\dots+c_nx_n=c_1g_1+\dots+c_ng_n$, so $g_1,\dots,g_n\in \dcl(g_1c_1+\dots+g_nc_n,C)$. If $\{g_1,\dots,g_n\}$ is a $C$-independent set, we get  $\uth(\tp(c_1g_1+\dots+c_ng_n/C))= n$. 
\end{proof}

\begin{prop} \label{not in scl implies U-rank omega}
Let $a\not \in \scl(\emptyset)$, then $\uth(\tp(a))=\omega$.
\end{prop}

\begin{proof}
By Theorem \ref{Main Theorem 3} (and Fact \ref{U-rank is smaller}), $\uth(\tp(a))\leq\omega$. 

Now we will show that $\tp(a/\emptyset)$ has forking extensions of $\uth$-rank $n$ for every $n$. Let $C= \{c_1,\dots,c_n\}$ and without loss of generality assume that $C$ is \thorn-independent from $a$. Let $g_1,\dots,g_n\in G$ and assume that $\{g_1,\dots,g_n\}$ is a $C\cup\{a\}$-independent set. Let $b=a+c_1g_1+\dots+c_ng_n$. Then $a$, $b\not \in \scl(\{c_1,\dots,c_n\})$. Thus $\uth(\tp(c_1g_1+\dots+c_ng_n/C\cup\{b\}))=n$ and since $a$ and $c_1g_1+\dots+c_ng_n$ are interdefinable over $b$, $\uth(\tp(a/C\cup \{b\}))=n$. Thus $\uth(\tp(a))=\omega$.
\end{proof}

\begin{cor}\label{being small causes thorn forking version 1}
If $\mathfrak R|_{\mathscr L}$ has a definable field structure and $a\in\scl(B)\setminus\scl(C)$, then $a\nthind_C B$.
\end{cor}

\begin{proof}
 We may assume $C=\emptyset$, as our hypotheses remain true after adding parameters to the language.  Since $a\in \scl(B)$, some formula in $\tp(a/B)$ defines a $G$-bound set, and  Lemma \ref{G-bound gives finite U-rank} implies that $\uth(a/B)$ is finite.  On the other hand, $\uth(a)=\omega$ by Lemma \ref{not in scl implies U-rank omega}.
\end{proof}

\subsection{Pairs of groups with no additional structure}
Assume now that $\mathscr L=\{+,0,<\}$. Thus $\mathfrak R|_{\mathscr L}$ is a divisible ordered abelian group.  Furthermore suppose that $G$ a subgroup of $R$.

\begin{defn}
Let $n>0$ and let $G/n=\{r\in R: nr\in G\}$.
\end{defn}

\begin{lem}
The group $G/n$ has \th-rank one.
\end{lem}

\begin{proof}
Recall that $G$ has \th-rank one. As $R$ is divisible and torsion-free, multiplication by $n$ is a definable bijection between $G/n$ and $G$, and thus the \th-rank of $G/n$ is one. 

\end{proof}

\begin{prop}\label{cosets enough}
$a\in\scl(B)$ if and only if there is $b\in\dcl(B)$ and $n\in\mathbb N^{>0}$ such that $a\in b+G/n$.
\end{prop}

\begin{proof}
Right to left is clear.

Now assume that $a\in\scl(B)$.  By Proposition \ref{basic small sets rock}, $a$ is contained in $\mathbb S$, a basic small set defined over $B$.  Let $\exists\vec y(\mathcal G(\vec y)\wedge\varphi(x,\vec y))$ be a formula defining $\mathbb S$. For each $\vec g$, $\varphi(R,\vec g)$ is a finite union of points and intervals.  However, if for any $\vec g$ in $G^k$, $\varphi(R,\vec g)$ contains a non-empty open interval, then $\mathbb S$ is not small.  Thus, we may reduce to the case where $\varphi(x,\vec y)$ is $x=f(\vec y)$, where $$f(\vec y)=b+\sum_{i=1}^k\frac{m_i}{n_i}y_i$$ for some $b\in\dcl(B)$, $m_i\in \mathbb Z$ and $n_i\in \mathbb N$. Let $n$ be the least common multiple of the $n_i$.  Thus $f(G^k)$ is contained in $b+G/n$, and $a\in f(G^k)$. 
\end{proof}

\begin{prop}
Let $a\in R$ be such that $a\not \in \scl(\emptyset)$. Then $\uth(\tp(a))= 2$. 
\end{prop}

\begin{proof}
By Proposition \ref{cosets enough}, every small subset of $R$ has \th-rank at most one, and by Corollary \ref{thorn dividing implies small}, a \th-forking extension of $\tp(a)$ must include a formula defining a small set.  Thus
$\uth(\tp(a))\leq 2$. It is easy to see that for $g\in G$, with $\tp(g)$ non-algebraic, and $g\thind a$, we get $\uth(\tp(a))=\uth(\tp(a/g))=\uth(\tp(a+g/g))$.  Now we claim that $a+g\thind g$.  If not, by Corollary \ref{thorn dividing implies small} we would have $a+g\in\scl(g)=\scl(\emptyset)$, and thus $a+g\in c+ G/n$ for some $c\in \dcl(\emptyset)$, by Proposition \ref{cosets enough}.  But then $a+g$, and hence $a$, would be in $\scl(\emptyset)$, a contradiction.  Thus $\uth(a)=\uth(a+g/g)=\uth(a+g)$, and it suffices to show that $\uth(a+g)=2$.

Consider the chain $\tp(a+g/\emptyset)\subset \tp(a+g/a)\subset \tp(a+g/a,g)$.  If we show that this is a \th-forking chain we will have shown that $\uth(a+g)\geq 2$, and thus equal to $2$. First note that $\tp(a+g/a)$ contains a formula saying $x\in G+a$.  This formula is true of $a+g$ and \th-divides over the empty set. Thus, $\tp(a+g/a)$ is a \th-forking extension of $\tp(a+g)$.

Second, note that $\tp(a+g/a,g)$ is algebraic, and hence to show that it is a \th-forking extension of $\tp(a+g/a)$, it suffices to show that the latter type is not algebraic.  But we chose $g\thind a$.  Thus $\tp(g/a)$ is not algebraic, and neither is $\tp(a+g/a)$.

\end{proof}

Now we get a corollary analogous to Corollary \ref{being small causes thorn forking version 1}:

\begin{cor}\label{being small causes thorn forking version 2}
If $\mathfrak R|_{\mathscr L}$ is an ordered group with no additional structure, and $a\in\scl(B)\setminus\scl(C)$, then $a\nthind_C B$.
\end{cor}
\begin{proof}
 By the previous proposition (after adding $C$ to the language), we see that $\uth(a/C)=2$.  On the other hand, by Proposition \ref{cosets enough}, we see that $a$ belongs to a set of \th-rank one defined over $B$, namely a coset of $G/n$ for some $n$.  Thus $\uth(a/B)$ is either zero or one.
\end{proof}

\begin{rem}
Note that we have shown that \th-forking in one variable is caused by falling into some coset of $G/n$ for some $n$.  This may be seen as an analogue of the fact from stable theories that the beautiful pair associated to a one-based theory is again one-based.
\end{rem}

\subsection{Small closure is a pregeometry}

\begin{cor}
If $\mathfrak R|_\mathscr L$ either is an ordered group with no additional structure or has a definable field structure, then the closure operator $\scl:\mathscr P(R)\to \mathscr P(R)$ defines a pregeometry.
\end{cor}

\begin{proof}
Let $C\subset R$ and let $a,b\in R$ be such that $a\in \scl(C\cup\{b\})\setminus \scl(C)$. Then  $\tp(a/C\cup\{b\})$ \th-forks over $C$ by either Corollary \ref{being small causes thorn forking version 1} or \ref{being small causes thorn forking version 2}. By symmetry, $\tp(b/C\cup\{a\})$ also \th-forks over $C$, so by Corollary \ref{thorn dividing implies small}, $b\in\scl(Ca)$.
\end{proof}

\end{document}